% Covering spaces of 3-orbifolds

\input amssym.def
\input amssym
\magnification=1100
\baselineskip = 0.25truein
\lineskiplimit = 0.01truein
\lineskip = 0.01truein
\vsize = 8.5truein
\voffset = 0.2truein
\parskip = 0.10truein
\parindent = 0.3truein
\settabs 12 \columns
\hsize = 5.4truein
\hoffset = 0.4truein
\font\ninerm=cmr9

\setbox\strutbox=\hbox{%
\vrule height .708\baselineskip
depth .292\baselineskip
width 0pt}
\font\caps=cmcsc10
\font\bigtenrm=cmr10 at 14pt

\font\smallmaths=cmmi9
\font\normalmaths=cmmi10

\def\sqr#1#2{{\vcenter{\vbox{\hrule height.#2pt
\hbox{\vrule width.#2pt height#1pt \kern#1pt
\vrule width.#2pt}
\hrule height.#2pt}}}}
\def\square{\mathchoice\sqr46\sqr46\sqr{3.1}6\sqr{2.3}4}

\centerline{\bigtenrm COVERING SPACES OF 3-ORBIFOLDS}
\tenrm
\vskip 14pt
\centerline{MARC LACKENBY}
\vskip 18pt

\vskip 12pt
\centerline{\caps Abstract}
\vskip 6pt

\ninerm
\textfont1 = \smallmaths
\baselineskip = 0.22truein

Let $O$ be a compact orientable 3-orbifold with non-empty singular locus 
and a finite volume hyperbolic structure. (Equivalently, the interior of $O$ is the quotient 
of hyperbolic 3-space by a lattice in ${\rm PSL}(2,{\Bbb C})$ with torsion.) Then we prove 
that $O$ has a tower of finite-sheeted covers $\{O_i\}$ with linear growth of 
$p$-homology, for some prime $p$. This means that the dimension of the first 
homology, with mod $p$ coefficients, of the fundamental group of $O_i$ grows 
linearly in the covering degree. The proof combines techniques from 3-manifold 
theory with group-theoretic methods, including the Golod-Shafarevich inequality 
and results about $p$-adic analytic pro-$p$ groups.
This has several consequences. Firstly, the fundamental group of $O$ has at 
least exponential subgroup growth. Secondly, the covers $\{O_i\}$ have positive
 Heegaard gradient. Thirdly, we use it to show that a group-theoretic 
conjecture of Lubotzky and Zelmanov would imply that $O$ has large fundamental 
group. This implication uses a new theorem of the author, which will appear 
in a forthcoming paper. These results all provide strong evidence for the 
conjecture that any closed orientable hyperbolic 3-orbifold with non-empty 
singular locus has large fundamental group.
Many of the above results apply also to 3-manifolds commensurable with an 
orientable finite-volume hyperbolic 3-orbifold with non-empty singular locus. 
This includes all closed orientable hyperbolic 3-manifolds with rank two 
fundamental group, and all arithmetic 3-manifolds.

\tenrm
\textfont1 = \normalmaths
\baselineskip = 0.25truein
\vskip 18pt

\footnote{}{MSC (2000): 57N10; 30F40; 20E07}

\tenrm
\centerline{\caps 1. Introduction}
\vskip 6pt

A central topic in 3-manifold theory is concerned
with a manifold's finite-sheeted covering spaces. The majority
of effort in this field was initially focused on
3-manifolds that are finitely covered by manifolds
that are well understood, in the hope that this would
provide illuminating information about the original
manifold. There are a number of excellent theorems
in this direction, for example, the result due
to Gabai, Meyerhoff and N. Thurston [5] that virtually
hyperbolic 3-manifolds are hyperbolic. Attention has
now turned to proving that hyperbolic 3-manifolds
always have finite covers with `nice' properties.
The central unresolved question in this direction
is the Virtually Haken Conjecture which proposes that any 
closed hyperbolic 3-manifold should be finitely covered by a Haken
3-manifold. There are even stronger versions of
this conjecture, but they all remain wide open at present.

There are other questions about finite-sheeted covering
spaces that arise naturally. A key one is:
how many finite-sheeted covering spaces does a 3-manifold
have, as a function of the covering degree? The answer
is not known in general, not even asymptotically. 
This question fits naturally into an emerging
area of group theory, which deals with `subgroup growth'.
This addresses the possible growth rates for the
number of finite index subgroups of a group, as a function
of their index. Linear groups play a central r\^ole
in this theory, particularly lattices in Lie groups.
However, discrete subgroups of ${\rm PSL}(2,{\Bbb C})$ or equivalently
the fundamental groups of orientable hyperbolic 3-orbifolds, remain
poorly understood.

There are other natural questions that can be asked,
which focus on the asymptotic behaviour of various
properties of the finite-sheeted covering
spaces. For example, one might wish to examine how
the following quantities can grow, as a function of
the covering degree: the rank of the fundamental group,
the rank and order of the first homology, the Heegaard genus,
the Cheeger constant and the first eigenvalue of the
Laplacian.
Indeed, a good understanding of these quantities is likely
to lead to progress on the Virtually Haken Conjecture,
by work of the author in [6].

The main tenet of this paper is that 3-orbifolds with non-empty
singular locus form a more tractable class than 3-manifolds,
when one is considering finite-sheeted covering spaces.
In particular, for discrete subgroups of ${\rm PSL}(2, {\Bbb C})$, it is
those with torsion that we can analyse most successfully.

Our first and central result is the following. It deals with
$d_p(\ \cdot \ )$. This is defined to be the dimension of
$H_1(\ \cdot \ ; {\Bbb F}_p)$, where ${\Bbb F}_p$ is the
field of prime $p$ order, and $\cdot$ is either a
group or a topological space. We say that a 
collection $\{ G_i \}$ of finite-index subgroups
of a finitely generated group $G$ has {\sl linear
growth of mod $p$ homology} if $\inf_i d_p(G_i)/[G:G_i] > 0$.

\vfill\eject
\noindent {\bf Theorem 1.1.} {\sl Let $O$ be a compact orientable
3-orbifold with non-empty singular locus and a finite-volume
hyperbolic structure. Then $O$
has a tower of finite-sheeted covers 
$\dots \rightarrow O_2 \rightarrow O_1 \rightarrow O$ where
$\{ \pi_1(O_i) \}$ has linear growth of mod $p$ homology,
for some prime $p$. Furthermore, one can ensure that
the following properties also hold:
\item{(i)} One can find such a sequence
where each $O_i$ is a manifold, and (when $O$ is closed)
another such sequence where each $O_i$ has non-empty singular locus.
\item{(ii)} Successive covers $O_{i+1} \rightarrow O_i$ are regular 
and have degree $p$.
\item{(iii)} For infinitely many $i$, $O_i \rightarrow O_1$
is regular.
\item{(iv)} One can choose $p$ to be any prime that divides the order
of an element of $\pi_1(O)$.

}

A slightly stronger version of the main part of this
theorem is as follows. 

\noindent {\bf Theorem 1.2.} {\sl Any finitely generated,
discrete, non-elementary subgroup of ${\rm PSL}(2, {\Bbb C})$
with torsion has a nested sequence of finite index subgroups
with linear growth of mod $p$ homology for some prime $p$.}

Linear growth of mod $p$ homology is a strong conclusion,
with several interesting consequences. For example,
one can use it to find good lower bounds
on the subgroup growth of the fundamental group of
a 3-orbifold. The {\sl subgroup growth function} $s_n(G)$,
for a finitely generated group $G$, is defined to be
the number of subgroups with index at most $n$. 
Of course, when $O$ is an orbifold, $s_n(\pi_1(O))$ simply
counts the number of covering spaces of $O$ (with given
basepoint) with degree at most $n$. It is
said to have {\sl polynomial growth} if there is
some constant $c$ such that $s_n(G) \leq n^c$ for each $n$.
It has {\sl at least exponential growth} if there is some
constant $c > 1$ such that $s_n(G) \geq c^n$ for each $n$.
Note that this is quite a strong form of exponential
growth: many authors just insist that $s_n(G) \geq c^n$
for infinitely many $n$.

Recall that an orbifold $O$ is {\sl geometric} if its interior
is a quotient $X/G$, where $X$ is a complete simply-connected
homogeneous Riemannian manifold and $G$ is a discrete
group of isometries. There are 8 possible model geometries
$X$ for compact 3-orbifolds. Furthermore, when $O$ is closed, it admits a geometry
modelled on at most one such model space. Thurston's
Orbifold Theorem ([1], [3]) states that a compact orientable
3-orbifold with non-empty singular locus and no bad
2-suborbifolds admits a `canonical decomposition into
geometric pieces'. Thus, it is natural to consider
geometric 3-orbifolds. In the following result,
we consider the subgroup growth of their fundamental
groups.

\noindent {\bf Theorem 1.3.} {\sl Let $O$ be a compact orientable 
geometric 3-orbifold
with non-empty singular locus. Then, the subgroup growth of $\pi_1(O)$ is
$$\cases{
\hbox{polynomial}, & if $O$ admits an $S^3$, ${\Bbb E}^3$, $S^2 \times {\Bbb E}$,
Nil or Sol geometry; \cr
\hbox{at least exponential}, & otherwise.}$$}

In fact, our techniques provide quite precise information about the
number of subnormal subgroups of $\pi_1(O)$. Recall that a 
subgroup $K$ of a group $G$ is {\sl subnormal}, denoted
$K \triangleleft \triangleleft \, G$, if there is a finite
sequence of subgroups $G = G_1 \triangleright G_2 \triangleright
\dots \triangleright G_n = K$, where each $G_i$ is normal
in $G_{i-1}$. Such groups arise naturally when one
considers towers of regular covers. For a finitely generated
group $G$, the number of subnormal subgroups of $G$
with index at most $n$ is denoted $s_n^{\triangleleft \triangleleft}
(G)$. It is known that $s_n^{\triangleleft \triangleleft}
(G)$ always grows at most exponentially (Theorem 2.3 of [16]). That is, there
is a constant $c$ depending on $G$ such that $s_n^{\triangleleft \triangleleft}
(G) \leq c^n$ for each $n$. We will show in Theorem 5.2
that, under the hypotheses of Theorem 1.3, $\pi_1(O)$ has
either polynomial or exponential subnormal subgroup
growth, depending on whether or not $O$ admits
an $S^3$, ${\Bbb E}^3$, $S^2 \times {\Bbb E}$,
Nil or Sol geometry.
In the hyperbolic case, this is rephrased as follows.

\noindent {\bf Theorem 1.4.} {\sl Any finitely generated, discrete,
non-elementary subgroup of ${\rm PSL}(2, {\Bbb C})$
with torsion has exponential subnormal subgroup growth.}

Theorem 1.1 also provides some new information about the behaviour
of the Heegaard Euler characteristic $\chi_-^h(O_i)$ of the manifold 
covering spaces $O_i$. 
Recall [6] that this is defined to be the negative of the
largest possible Euler characteristic for a Heegaard
surface for $O_i$. The
{\sl Heegaard gradient} of a sequence $\{ O_i \rightarrow O \}$
of finite-sheeted (manifold) covers is defined to be
$\inf_i \chi_-^h(O_i)/{\rm degree}(O_i \rightarrow O)$.
Since $\chi_-^h(O_i)$ is bounded below by a linear function of  
$d_p(O_i)$, for any prime $p$, linear growth of mod $p$ homology
implies positive Heegaard gradient. Thus, Theorem 1.1
has the following immediate consequence.

\noindent {\bf Corollary 1.5.} {\sl Let $O$ be a compact
orientable 3-orbifold with non-empty singular locus and
a finite-volume hyperbolic structure. Then $O$
has a nested sequence of finite-sheeted manifold covers with
positive Heegaard gradient.}

This is related to a proposed approach to the Virtually Haken
Conjecture, which arises from the following theorem of the author [6].
(The theorem below is not stated precisely this way in [6], but
this version is an immediate consequence of Theorem 1.4 of [6].)

\noindent {\bf Theorem 1.6.} {\sl Let $M$ be a closed orientable irreducible
3-manifold, and let $\{ M_i \rightarrow M \}$
be a nested sequence of finite-sheeted regular covering spaces of $M$.
Suppose that
\item{(i)} the Heegaard gradient of $\{ M_i \rightarrow M \}$ 
is positive; and
\item{(ii)} $\pi_1(M)$ does not have Property $(\tau)$ with
respect to $\{ \pi_1(M_i) \}$.

\noindent Then $M_i$ is Haken for all sufficiently large $i$.}

Property $(\tau)$, referred to in the above theorem, is an important 
concept from group theory, with links to many areas of mathematics, including
representation theory, graph theory and differential geometry [12].
Lubotzky and Sarnak [13] conjectured that a closed hyperbolic 3-manifold
$M$ should always have a sequence of finite-sheeted
covering spaces where (ii) holds. A key question is
whether this is true for the covers in Corollary 1.5.
A positive answer would settle the Virtually Haken Conjecture
for compact orientable hyperbolic 3-orbifolds with non-empty singular locus
(by setting $M$ in Theorem 1.6 to be the first manifold cover $O_1$). 
Indeed, the following recent theorem of the author
would provide a much stronger conclusion. Recall that a
group is {\sl large} if it has a finite index subgroup
that admits a surjective homomorphism onto a non-abelian
free group.

\noindent {\bf Theorem 1.7.} {\sl Let $G$ be a finitely presented group,
let $p$ be a prime and suppose that $G \geq G_1 \triangleright G_2 \triangleright \dots$
is a nested sequence of finite index subgroups, such
that $G_{i+1}$ is normal in $G_{i}$ and has index a power of $p$, for each $i$. 
Suppose that $\{ G_i \}$ has linear growth of mod $p$ homology. Then,
at least one of the following must hold:
\item{(i)} $G$ is large;
\item{(ii)} $G$ has Property $(\tau)$ with respect to $\{G_i \}$.

}

This is, in fact, a slightly abbreviated and slightly
weaker version of the main theorem (Theorem 1.1) of [8].
Set $G$ to be $\pi_1(O)$ and let $G_i$ be $\pi_1(O_i)$,
where $\{ O_i \rightarrow O\}$ is one of the sequences of
covering spaces in Theorem 1.1.
We know from Theorem 1.1 that $\{ G_i \}$ has linear growth
of mod $p$-homology.
Thus, the central question is: does $G$ have Property
$(\tau)$ with respect to $\{ G_i \}$? We conjecture
that we may pick $\{ G_i \}$ so that it does not. In fact, we will see that this
would follow from the following recent conjecture
of Lubotzky and Zelmanov [14], which we have termed the
GS-$\tau$ Conjecture.

\noindent {\bf Conjecture 1.8.} (GS-$\tau$ Conjecture) {\sl Let $G$ be  
a group with finite presentation $\langle X | R \rangle$,
and let $p$ be a prime.
Suppose that $d_p(G)^2/4 > |R| - |X| + d_p(G)$.
Then $G$ does not have Property $(\tau)$
with respect to some infinite nested sequence $\{ G_i \}$
of normal subgroups with index a power of $p$.}

The point is that the pro-$p$ completion of $G$ has
a pro-$p$ presentation with $d_p(G)$ generators and $(|R| - |X| + d_p(G))$
relations. Thus, the condition in the GS-$\tau$ Conjecture asserts
that the Golod-Shafarevich inequality does not hold
for this pro-$p$ presentation. The theory of such pro-$p$
groups is advanced and there is some hope that
it may be applied and developed to prove this conjecture.
Indeed, the Golod-Shafarevich inequality will play
a central r\^ole in this paper. 

We will show that the following is a consequence of Theorems 1.1 and 1.7.

\noindent {\bf Theorem 1.9.} {\sl The GS-$\tau$
Conjecture implies that the fundamental group of 
every closed hyperbolic 3-orbifold with non-empty singular locus 
is large.}

It therefore seems that the evidence for the following
conjecture is stacking up.

\noindent {\bf Conjecture 1.10.} {\sl The fundamental group
of any closed hyperbolic 3-orbifold with non-empty singular locus
is large.}

Indeed, Theorems 1.1 and 1.3 are already pointing in this direction.
This is because linear growth of mod $p$ homology for some nested
sequence of finite-index subgroups and (at least)
exponential subgroup growth are both strong properties
of a group, which are enjoyed when the group is
large.

It is worth pointing out that the GS-$\tau$ Conjecture is not the only possible approach to
Conjecture 1.10. Another possible route, explained in [8],
is via error-correcting codes.

Another piece of evidence for Conjecture 1.10 comes
from work of the author in [7]. There, a sequence
of results about 3-orbifolds is established, one of
which is the following.

\noindent {\bf Theorem 1.11.} (Theorem 3.6 of [7]) 
{\sl Let $O$ be a compact orientable 3-orbifold (with possibly
empty singular locus), and
let $K$ be a knot in $O$, disjoint from its singular
locus, such that $O - K$
admits a finite-volume hyperbolic structure. 
For any integer $n$, let $O(K,n)$ denote the
3-orbifold obtained from $O$ by adjoining
a singular component along $K$ with order $n$.
Then, for infinitely many
values of $n$, $\pi_1(O(K,n))$ is large.}

Of course, Conjecture 1.10 is weaker than the old conjecture that 
the fundamental
group of any closed hyperbolic 3-manifold is large, because any closed
hyperbolic 3-orbifold is finitely covered by a hyperbolic 3-manifold,
by Selberg's Lemma. But the main purpose
of this paper is to demonstrate that covering spaces of
3-orbifolds with non-empty singular locus
are more tractable than the manifold case,
and so there is some chance that Conjecture 1.10 may be
more likely to be true.

The slightly stronger version of Theorem 1.7 that appears
in [8] has the following consequence. (See [8] for an explanation of this deduction.)

\noindent {\bf Theorem 1.12.} {\sl Let $G$ be a finitely
presented group, and let $p$ be a prime. Suppose that
$G$ has an infinite nested sequence of subnormal subgroups
$\{ G_i \}$, each with index a power of $p$
and with linear growth of mod $p$ homology.
Then $G$ has such a sequence that also has Property $(\tau)$.}

Combining this with Theorem 1.1, we have the following
corollary.

\noindent {\bf Theorem 1.13.} {\sl Let $O$ be a compact
orientable 3-orbifold with non-empty singular locus and
a finite-volume hyperbolic structure. Then $O$ has
a nested sequence of finite-sheeted
covers with Property $(\tau)$.}

It is worth pointing out that the above results
remain true for hyperbolic 3-manifolds that are commensurable
with hyperbolic 3-orbifolds with non-empty singular locus.

\noindent {\bf Theorem 1.14.} {\sl Let $M$ be a 3-manifold
that is commensurable with a compact orientable finite-volume hyperbolic 3-orbifold $O$
with non-empty singular locus. Let $p$ be a prime that
divides the order of a torsion element of $\pi_1(O)$. Then $\pi_1(M)$ has
exponential subnormal subgroup growth, and $M$ has a nested
sequence of finite-sheeted covers that have linear growth
of mod $p$ homology and have Property $(\tau)$.}

There are many examples of such 3-manifolds $M$. A large class of examples is the
orientable hyperbolic 3-manifolds with rank two fundamental group,
which are 2-fold regular
covers of hyperbolic 3-orbifolds with non-empty singular
locus (Corollary 5.4.2 of [19]).

\noindent {\bf Corollary 1.15.} {\sl If $M$ a closed
orientable hyperbolic 3-manifold with rank two
fundamental group, then the conclusions of Theorems 1.14
apply to $M$, for $p=2$.}

Another important class of 3-manifolds to which one may apply
these results are arithmetic hyperbolic 3-manifolds.
The following theorem appears in [9].

\noindent {\bf Theorem 1.16.} {\sl Let $M$ be 
an arithmetic hyperbolic 3-manifold. Then
$M$ is commensurable with an arithmetic hyperbolic 3-orbifold,
with fundamental group that contains
${\Bbb Z}/2{\Bbb Z}\times{\Bbb Z}/2{\Bbb Z}$.}

Thus, we have the following result.

\noindent {\bf Corollary 1.17.} {\sl If $M$ is an arithmetic
hyperbolic 3-manifold, then the conclusions of Theorem 1.14 apply
to $M$ for $p =2$.}

In this paper, we assume some familiarity with the basic
theory of orbifolds. In particular, we take as given the
following terminology: singular locus, covering
space, fundamental group. For an explanation of these terms
and introduction to orbifolds, we suggest [3] and [18]
as helpful references.

\vskip 18pt
\centerline{\caps 2. Covering spaces without vertices}
\vskip 6pt

The main focus of this paper will be 3-orbifolds
that admit a hyperbolic structure and that have non-empty singular locus. 
It turns out that it is easiest to deal with 3-orbifolds that
have singular locus containing no vertices. Our goal in this section
is to show that one can always arrange this to be the
case, by first passing to a finite-sheeted cover.

\noindent {\bf Proposition 2.1.} {\sl Let $O$ be a compact orientable hyperbolic
3-orbifold with non-empty singular locus. Then $O$ has a 
finite-sheeted cover $\tilde O$ with singular locus that is non-empty and that contains no vertices.
Moreover, for any prime $p$ that divides the order of an element of $\pi_1(O)$,
we may arrange that the order of each component of this
singular locus is $p$. In addition, we may ensure that
there is a degree $p$ regular cover $M \rightarrow \tilde O$ 
where $M$ is a manifold, and where the composite cover
$M \rightarrow \tilde O \rightarrow O$ is regular.}

\noindent {\sl Proof.} Since $O$ is hyperbolic, its fundamental
group is realised as a subgroup $G$ of ${\rm PSL}(2, {\Bbb C})$.
By Selberg's Lemma, $G$ has a finite index normal subgroup $K$ that
is torsion free. Let $M$ be the regular covering space of $O$ corresponding
to $K$; this is a manifold. Let $\mu$ be an element of $G$ that has finite
order. We may choose $\mu$ so that it has order $p$.
We claim that the covering space $\tilde O$ of $O$ corresponding 
to $K \langle \mu \rangle$ has non-empty singular locus that contains 
no vertices. Since $K$ is normal in $G$, it is also normal in
$K \langle \mu \rangle$. So, $M \rightarrow \tilde O$ is a regular
cover, with covering group $K \langle \mu \rangle / K$, which is
cyclic of order $p$. Hence, $\tilde O$ is the quotient of $M$ by a finite
order orientation-preserving diffeomorphism. The singular locus of $\tilde O$ is
therefore a 1-manifold: it has no vertices. It is non-empty, because
the fundamental group of $\tilde O$, namely $K \langle \mu \rangle$,
contains $\mu$ which has finite order. Since the order of the
diffeomorphism is $p$, which is prime, the order of every component of the
singular locus is $p$. $\square$

We now introduce some convenient terminology.
If $p$ is a prime and $O$ is a 3-orbifold, then
we let ${\rm sing}^c_p(O)$ denote those simple
closed curve components of the singular locus
of $O$ with singularity order that is a multiple
of $p$. 

\vskip 18pt
\centerline{\caps 3. Covering spaces with large mod $p$ homology}
\vskip 6pt

The goal of this section is to prove the
following result.

\noindent {\bf Theorem 3.1.} {\sl Let $O$ be
an orientable finite-volume
hyperbolic 3-orbifold. Let $p$ be a prime such that
${\rm sing}^c_p(O)$ is non-empty. Then there is a
finite-sheeted cover $\tilde O \rightarrow O$,
such that $d_p(\pi_1(\tilde O)) \geq 11$ 
and where ${\rm sing}^c_p(\tilde O)$ is non-empty.
In addition, we may ensure that there is
a regular covering space $\tilde O' \rightarrow \tilde O$
which has degree $p$ or $1$ and where the
composite cover $\tilde O' \rightarrow O$ is
regular.}

The significance of the number 11 in the above theorem
will be made clear in Theorem 4.1 and its proof.
The proof of Theorem 3.1 relies on the following
result, which is due to Lubotzky, and is of independent interest.

\noindent {\bf Theorem 3.2.} {\sl Let $O$ be an orientable finite-volume
hyperbolic 3-orbifold (with possibly empty singular locus).
Then, for any prime $p$,
$$\sup \{ d_p(K) : K \hbox{ is a finite index normal subgroup of }
 \pi_1(O) \} = \infty.$$}

The proof of this requires three theorems, which we quote.
The first is a consequence
of Nori-Weisfeiler's Strong Approximation Theorem
and the Lubotzky Alternative
(see Corollary 18 of Window 9 in [16] for example.)

\noindent {\bf Theorem 3.3.} {\sl Let $G$
be a finitely generated linear group that
is not virtually soluble.
Then, for any prime $p$, 
$$\sup \{ d_p(K) : K \hbox{ is a finite index subgroup of }
 G \} = \infty.$$}

The subgroups provided by this theorem need not be
normal. In fact, there exist groups $G$ satisfying the
hypotheses of Theorem 3.3 but where
$$\sup \{ d_p(K): K \hbox{ is a finite index normal subgroup of }
 G \}$$
is finite, for infinitely many primes $p$. An example is
${\rm SL}(n, {\Bbb Z})$, for any $n \geq 3$ and any prime $p$
not dividing $n$ (see Proposition 1.4 of [11]).
So, more work is necessary before Theorem 3.2
can be proved. The second result we quote can be found in [15].

\noindent {\bf Theorem 3.4.} {\sl Let $G$ be
a finitely generated group, let $p$ be a prime
and let $\hat G_{(p)}$ be the pro-$p$ completion of $G$.
Then the following are equivalent:
\item{1.} $\hat G_{(p)}$ is $p$-adic analytic;
\item{2.} the supremum of $d_p(K)$, as $K$ ranges over all
characteristic subgroups of $G$ with index a power of $p$,
is finite;
\item{3.} the supremum of $d_p(K)$, as $K$ ranges over all
normal subgroups of $G$ with index a power of $p$,
is finite.

}

The relationship between $p$-adic analytic pro-$p$ groups
and 3-manifolds can be seen in the following result, which
is due to Lubotzky [10].

\noindent {\bf Theorem 3.5.} {\sl Let $M$ be a
compact orientable 3-manifold and let $p$ be a prime.
If $d_p(M) \geq 4$, then the pro-$p$ completion of 
$\pi_1(M)$ is not $p$-adic analytic.}

The following is elementary and fairly well known.

\noindent {\bf Proposition 3.6.} {\sl Let $G$ be
a finitely generated group, and let $K$ be a normal
subgroup. Then, $d_p(G) \leq d_p(K) + d_p(G/K)$.}

\noindent {\sl Proof.} Let $G'$ be $[G,G]G^p$,
which is the subgroup generated by the commutators
and $p^{\rm th}$ powers of $G$.
Define $K'$ similarly. Consider the exact
sequence
$$1 \rightarrow {KG' \over G'} \rightarrow {G \over G'} \rightarrow
{G \over KG'} \rightarrow 1.$$
A set of elements in $G/G'$ that maps to a generating
set for $G/KG'$, together with a generating set
for $KG'/G'$, forms a generating set for $G/G'$.
Hence, writing $d( \ )$ for the minimal number of
generators for a group, 
$$d_p(G) = d(G/G') \leq d(KG'/G') + d(G/KG').$$
Now, $d(G/KG') = d_p(G/K)$. Also, $KG'/G'$ is
isomorphic to $K/(K \cap G')$, which is a quotient
of $K/K'$. Hence, $d(KG'/G') \leq d(K/K') = d_p(K)$.
The required inequality now follows. $\square$

\noindent {\sl Proof of Theorem 3.2.}
Let $O$ be a compact orientable 3-orbifold with a finite-volume
hyperbolic structure. Then $\pi_1(O)$ is realised as
a lattice in ${\rm PSL}(2, {\Bbb C})$. Therefore,
by Selberg's Lemma, $\pi_1(O)$ has a finite index subgroup
which is torsion free. This corresponds to a manifold
covering space of $O$. By Theorem 3.3, this has a finite-sheeted
cover $M$ where $d_p(M) \geq 4$, say. Let $M' \rightarrow O$
be the covering corresponding to the intersection of all
conjugates of $\pi_1(M)$ in $\pi_1(O)$. Then $M' \rightarrow M$
and $M' \rightarrow O$ are finite-sheeted regular covers.

By Theorem 3.5, the pro-$p$ completion of $\pi_1(M)$ is not $p$-adic analytic.
Hence, by Theorem 3.4, $M$ has finite-sheeted covers $M_i$
such that $\pi_1(M_i)$ is characteristic in $\pi_1(M)$
and has index a power of $p$, and where $d_p(M_i)$
tends to infinity. Let $M'_i$ be the cover of $M$
corresponding to the subgroup $\pi_1(M') \cap \pi_1(M_i)$.
This is the intersection of two normal subgroups of 
$\pi_1(M)$ and so is normal in $\pi_1(M)$. Hence, $M'_i$ regularly
covers $M$, $M'$ and $M_i$. Now, by Proposition
3.6, $d_p(M'_i) \geq d_p(M_i) - d_p(\pi_1(M_i)/ \pi_1(M'_i))$.
But
$\pi_1(M_i)/\pi_1(M'_i) = \pi_1(M_i) / (\pi_1(M_i) \cap \pi_1(M'))
= \pi_1(M_i) \pi_1(M') / \pi_1(M')$, which is a subgroup
of $\pi_1(M)/ \pi_1(M')$. There are only finitely many such
subgroups and hence their $d_p$ is uniformly bounded above.
So, $d_p(M'_i)$ tends to infinity.

Now, $\pi_1(M') / \pi_1(M'_i) = \pi_1(M') / (\pi_1(M_i) \cap \pi_1(M'))
= \pi_1(M') \pi_1(M_i) / \pi_1(M_i)$, which is a subgroup
of $\pi_1(M)/\pi_1(M_i)$. This is a finite $p$-group, and hence so 
is $\pi_1(M') / \pi_1(M'_i)$. Thus, $\pi_1(M'_i)$ corresponds to a
finite index normal subgroup of the pro-$p$ completion of
$\pi_1(M')$. Since $d_p(M'_i)$ tends to infinity, this
pro-$p$ completion is not $p$-adic analytic, by Theorem 3.4,
and so, by Theorem 3.4, $\pi_1(M')$ has a sequence of
characteristic subgroups $K_i$, each with index a power of $p$,
and where $d_p(K_i)$ tends to infinity. Since these are
characteristic in $\pi_1(M')$, which is normal in $\pi_1(O)$,
they are therefore normal in $\pi_1(O)$. These are the
required subgroups. $\square$

Note that, in the proof of Theorem 3.2, we only used once the
hypothesis that $\pi_1(O)$ is the fundamental group of an orientable
finite-volume hyperbolic 3-orbifold. This was used to show that it
has a finite index subgroup that is not $p$-adic analytic. Once this
is known, the remainder of the proof is purely group-theoretic.

\noindent {\sl Proof of Theorem 3.1.}
Let $\mu$ be an element of $\pi_1(O)$ that
forms a meridian for a component $C$ of 
${\rm sing}^c_p(O)$. Then
$\mu^n$ has order $p$ for some positive integer $n$.
Let $K$ be a finite index normal subgroup of $\pi_1(O)$,
as in Theorem 3.2, where $d_p(K)$ is large (more
than $10p$ will suffice). Let $\tilde O'$ be
the regular covering space of $O$ corresponding to this
subgroup.

If $\mu^n$ lies in $K$, then we set $\tilde O$
to be $\tilde O'$. Note that ${\rm sing}^c_p(\tilde O)$
is non-empty because the inverse image of
$C$ in $\tilde O$ contains at least one singular
component such that $p$ divides its singularity order, 
as $\mu^n$ lies in $K$.

So, suppose that $\mu^n$ does not lie in $K$.
We will consider the covering space
$\tilde O$ corresponding to the subgroup $K \langle \mu^n \rangle$.
Note that again ${\rm sing}^c_p(\tilde O)$ is
non-empty. Note also that $\tilde O' \rightarrow \tilde O$
has degree equal to the order of $K \langle \mu^n \rangle /K$,
which is $p$.
We want to show that, when $d_p(K)$ is
large, then so is $d_p(K\langle \mu^n \rangle)$.
In particular, when $d_p(K) \geq 10p+1$,
then $d_p(K\langle \mu^n \rangle) \geq 11$.
This is an immediate consequence of the 
following well known proposition (setting
$H = K \langle \mu^n \rangle$). $\square$

\noindent {\bf Proposition 3.7.} {\sl Let $H$ be a finitely
generated group, and let $K$ be a subnormal subgroup with index
a power of a prime $p$. Then
$$d_p(K) - 1 \leq [H:K] (d_p(H) - 1).$$}

\noindent {\sl Proof.} By a straightforward induction,
we may reduce to the case where $K$ is a normal
subgroup of $H$. Let $K'$ denote $[K,K]K^p$,
and define $H'$ similarly. Then $K/K'$ is 
isomorphic to $H_1(K; {\Bbb F}_p)$. Now, $K'$ is characteristic
in $K$, which is normal in $H$, and hence $K'$ is
normal in $H$. Thus, $H/K'$ is a finite $p$-group.
It is a well known fact that any finite $p$-group
has a generating set with size 
equal to the dimension of its first homology with
${\Bbb F}_p$-coefficients, in this case
$d_p(H/K')$. But $d_p(H/K') = d_p(H/K'H') = d_p(H)$.
Now, $K/K'$ is a subgroup of $H/K'$ with index
$[H:K]$. Applying the Reidemeister-Schreier
process to this subgroup, we obtain a generating
set for $K/K'$ with size $[H:K] (d_p(H) - 1) + 1$.
This is therefore an upper bound for $d_p(K)$. 
The required inequality follows.
$\square$

\vskip 18pt
\centerline{\caps 4. Linear growth of the number of singular components}
\vskip 6pt

We now come to the central result of this paper, from which
all later analysis of hyperbolic 3-orbifolds will follow.

\noindent {\bf Theorem 4.1.} {\sl Let $O$ be a compact orientable 3-orbifold,
with boundary a \break (possibly empty) union of tori, and 
with singular locus that is a link. Let
$p$ be a prime that divides the order of a component
$C$ of the singular locus, and let $\langle \! \langle \pi_1(\partial N(C))
\rangle \! \rangle$ be the
subgroup of $\pi_1(O)$ normally generated by 
$\pi_1(\partial N(C))$. Suppose that $d_p(\pi_1(O)) \geq 11$.
Then $\pi_1(O)$ has an infinite nested sequence of
finite index subgroups $\{ G_i \}$ such that
\item{(i)} each $G_i$ contains $\langle \! \langle \pi_1(\partial N(C))
\rangle \! \rangle$;
\item{(ii)} each $G_{i+1}$ is normal in $G_i$ and has index $p$; and
\item{(iii)} infinitely many $G_i$ are normal in
$\pi_1(O)$ and have index a power of $p$.

}

The significance of this theorem comes from the following
proposition.

\noindent {\bf Proposition 4.2.} {\sl Let $O$ be a 3-orbifold,
let $p$ be a prime, and let $C$ be some component of
${\rm sing}^c_p(O)$. Let $\tilde O \rightarrow O$ be a
finite-sheeted cover, such that $\pi_1(\tilde O)$
contains the subgroup of $\pi_1(O)$ normally
generated by $\pi_1(\partial N(C))$. Then
$|{\rm sing}^c_p(\tilde O)| \geq {\rm degree}(\tilde O \rightarrow O)$.}

\noindent {\sl Proof.}
The condition that $\pi_1(\tilde O)$ contains
the subgroup of $\pi_1(O)$ normally generated
by $\pi_1(\partial N(C))$
forces the inverse image of $\partial N(C)$
in $\tilde O$ to be a disjoint union of copies of 
$\partial N(C)$. Each copy bounds a component
of ${\rm sing}^c_p(\tilde O)$. $\square$

We will also need the following key fact,
which will be used both in the proof of
Theorem 4.1 and in its applications later.

\noindent {\bf Proposition 4.3.} {\sl Let $O$ be a compact
orientable 3-orbifold with singular locus consisting
of a link. Then, for any prime $p$, $d_p(\pi_1(O)) \geq |{\rm sing}^c_p(O)|$.}

\noindent {\sl Proof.} We may remove the components of
the singular locus of $O$ that have order coprime to $p$,
and replace them by manifold points. This is because
such an operation changes neither $d_p(\pi_1(O))$ nor
${\rm sing}^c_p(O)$. Thus, we may assume that every
component of the singular locus lies in ${\rm sing}^c_p(O)$.

Let $M$ be $O - {\rm int}(N({\rm sing}^c_p(O))$.
Then, it is a well-known consequence of Poincar\'e duality
that $d_p(M)$ is at least ${1 \over 2}d_p(\partial M)$,
which is at least $|{\rm sing}^c_p(O)|$. One obtains $\pi_1(O)$ from
$\pi_1(M)$ by adding relations that are each $p^{\rm th}$ powers
of words. This does not affect $d_p$. We therefore
obtain the required inequality. $\square$

We will need the following well known theorem. It deals
with groups $\Gamma$ whose pro-$p$ completion has
a (minimal) presentation that does not satisfy the Golod-Shafarevich
condition. A proof can be found in [10]. In fact, Theorem 1.2 of
[10] is a significantly stronger result.

\noindent {\bf Theorem 4.4.} {\sl Let $\langle X | R \rangle$ be a finite
presentation of a group $\Gamma$, and let $p$ be a prime. Suppose that 
$$d_p(\Gamma)^2/4 > |R| - |X| + d_p(\Gamma).$$
Then $\Gamma$ has an infinite nested sequence of subgroups
$\Gamma = \Gamma_1 \geq \Gamma_2 \geq \dots$,
where $\Gamma_{i+1}$ is normal in $\Gamma_i$ and has
index $p$. Furthermore, infinitely many $\Gamma_i$
are normal in $\Gamma$.}

\noindent {\bf Remark.} The normal subgroups of $\Gamma$ can, in fact, be taken
to be the $p$-lower central series, which is defined by
setting $\Gamma_1 = \Gamma$ and $\Gamma_{j+1} = [\Gamma_j,\Gamma]\Gamma_j^p$.
The quotients $\Gamma_j / \Gamma_{j+1}$ are elementary abelian
$p$-groups, and therefore we may interpolate between 
$\Gamma_j$ and $\Gamma_{j+1}$ by a sequence of subgroups,
where each is normal in its predecessor and has index $p$.

\noindent {\sl Proof of Theorem 4.1.} 
Let $O'$ be the orbifold with the same underlying
manifold as $O$ but with singular locus equalling
${\rm sing}^c_p(O)$. In other words, we remove
from $O$ those components of the singular locus
with order coprime to $p$, and replace them by manifold points. Then there is a
natural map $O \rightarrow O'$ which induces
a surjective homomorphism $\pi_1(O) \rightarrow \pi_1(O')$
and an isomorphism $H_1(\pi_1(O); {\Bbb F}_p) \rightarrow
H_1(\pi_1(O'); {\Bbb F}_p)$. In particular,
$d_p(\pi_1(O')) \geq 11$. If we can show that the
theorem holds for $O'$, then it also holds for
$O$. This is because the nested sequence of
subgroups of $\pi_1(O')$ pulls back to give
a corresponding sequence of subgroups of $\pi_1(O)$
with the required properties. Hence, by replacing
$O$ by $O'$, we may assume that $p$ divides
the order of every component of the singular
locus of $O$.

Removing the singular locus from $O$ gives
a compact orientable 3-manifold $M$ with boundary a 
non-empty union of tori. Its fundamental
group therefore has a presentation where the number
of relations equals the number of generators minus 1. 
Let $\Gamma$ be $\pi_1(O)/\langle \! \langle \pi_1(\partial N(C))
\rangle \! \rangle$.
Then, $\Gamma$ has a presentation $\langle X | R \rangle$
where $|R| - |X| = |{\rm sing}^c_p(O)|$, which is obtained 
from the presentation of $\pi_1(M)$ by
killing a pair of generators for $\pi_1(\partial N(C))$,
and adding a relation for each remaining component of
the singular locus of $O$. By Proposition 4.3,
$d_p(\pi_1(O)) \geq |{\rm sing}^c_p(O)|$, and hence
$|R| - |X| \leq d_p(\pi_1(O))$.
Now, $d_p(\Gamma) \geq d_p(\pi_1(O)) - 2$,
since $\Gamma$ is obtained from $\pi_1(O)$ by
adding two relations.
This implies that
$$\eqalign{
d_p(\Gamma)^2/4 - d_p(\Gamma) + |X| - |R| 
&\geq d_p(\Gamma)^2/4 - 2d_p(\Gamma) - 2  \cr
&\geq (d_p(\pi_1(O)) - 2)^2/4 - 2 (d_p(\pi_1(O)) - 2) - 2 \cr
& \geq 9^2/4 - 18 - 2 > 0.}$$
Note that the second and third inequalities hold because
$x^2/4 - 2x - 2$ is an increasing function of $x$
when $x \geq 4$. 
So, by Theorem 4.4, $\Gamma$ has an infinite nested
sequence of subgroups $\Gamma = \Gamma_1 \geq \Gamma_2 \geq \dots$, 
where each $\Gamma_{i+1}$ is normal in $\Gamma_i$ and has index $p$.
Furthermore, infinitely many $\Gamma_i$ are normal in
$\Gamma$.
These pull back to give the required subgroups of $\pi_1(O)$.
$\square$

\vskip 18pt
\centerline{\caps 5. Growth of homology and subgroup growth}
\vskip 6pt

We can now prove the main theorem stated in the
introduction.

\noindent {\bf Theorem 1.1.} {\sl Let $O$ be a compact orientable
3-orbifold with non-empty singular locus and a finite-volume
hyperbolic structure. Then $O$
has a tower of finite-sheeted covers 
$\dots \rightarrow O_2 \rightarrow O_1 \rightarrow O$ where
$\{ \pi_1(O_i) \}$ has linear growth of mod $p$ homology,
for some prime $p$. Furthermore, one can ensure that
the following properties also hold:
\item{(i)} One can find such a sequence
where each $O_i$ is a manifold, and (when $O$ is closed)
another such sequence where each $O_i$ has non-empty singular locus.
\item{(ii)} Successive covers $O_{i+1} \rightarrow O_i$ are regular 
and have degree $p$.
\item{(iii)} For infinitely many $i$, $O_i \rightarrow O_1$
is regular.
\item{(iv)} One can choose $p$ to be any prime that divides the order
of an element of $\pi_1(O)$.

}

\noindent {\sl Proof.} Consider first
the main case, where $O$ is closed. 
By Proposition 2.1, we may
pass to a finite cover $O'$ of $O$ with non-empty
singular locus, each component of which is a simple
closed curve with order $p$. Moreover, we may
ensure that $O'$ has a degree $p$ regular cover that
is a manifold $M$. By Theorem 3.1,
there is a finite-sheeted cover $O''$ of $O'$,
again where the singular locus is a non-empty
collection of simple closed curves with order $p$, 
but also where $d_p(\pi_1(O'')) \geq 11$.
Now apply Theorem 4.1 and Proposition 4.2 to deduce
the existence of a tower of finite covers $O_i \rightarrow O''$
such that $|{\rm sing}^c_p(O_i)| \geq {\rm degree}(O_i \rightarrow O'')$.
The covers $O_i \rightarrow O$ will be those required
by the theorem.
Proposition 4.3 gives that $d_p(\pi_1(O_i)) \geq
|{\rm sing}^c_p(O_i)|$. Thus, $\{ \pi_1(O_i) \}$ has
linear growth of mod $p$ homology, as required. 
By (ii) of Theorem 4.1, we may ensure that, for each
$i$, $O_{i+1} \rightarrow O_i$ is regular and has
degree $p$. By (iii) of Theorem 4.1, for
infinitely many $i$, $O_i \rightarrow O''$ is
regular, and hence so is $O_i \rightarrow O_1$.

We now show
how to find another such sequence of covers consisting of manifolds.
Now, $O_i$ has a manifold
cover $M_i$ corresponding to the subgroup $\pi_1(O_i) \cap
\pi_1(M)$. This is a regular cover, with degree
equal to $[\pi_1(O_i): \pi_1(O_i) \cap \pi_1(M)]
= [\pi_1(O_i) \pi_1(M) : \pi_1(M)] = p$,
since $\pi_1(O_i) \pi_1(M)/\pi_1(M)$ is a non-trivial
subgroup of $\pi_1(O')/\pi_1(M)$, which is cyclic
of order $p$. The fact that $\{ \pi_1(M_i) \}$
has linear growth of mod $p$ homology is a consequence of
Proposition 3.6 (letting $G = \pi_1(O_i)$ and
$K = \pi_1(M_i)$). Conclusions (ii) and (iii) of Theorem 1.1
apply to the sequence $\{ M_i \}$ because they
apply to $\{ O_i \}$.

Suppose now that $O$ has non-empty boundary. 
By Selberg's Lemma, we may pass to a finite-sheeted
regular manifold cover $M$. By a result of Cooper,
Long and Reid (Theorem 1.3 of [4], see also [7], [2] and [17]), 
$\pi_1(M)$ is large: it has a finite-index normal subgroup
that admits a surjective homomorphism onto
a free non-abelian group $F$. Let $O_1$ be the
covering space corresponding to this subgroup.
Now, $F$ has
a nested sequence of finite-index subgroups $\{ F_i \}$
with linear growth of mod $p$ homology.
We may ensure that each $F_{i+1}$ is normal in $F_i$
with index $p$, that $F_1 = F$, and that infinitely many $F_i$ are
normal in $F$.
The inverse image of these subgroups in
$\pi_1(O)$ correspond to the required covering
spaces $O_i$ of $O$. $\square$

A slight extension of Theorem 1.1 is the following.

\noindent {\bf Theorem 1.2.} {\sl Any finitely generated,
discrete, non-elementary subgroup of ${\rm PSL}(2, {\Bbb C})$
with torsion has a nested sequence of finite index subgroups
with linear growth of mod $p$ homology for some prime $p$.}

\noindent {\sl Proof.}
Let $O$ be the quotient of ${\Bbb H}^3$ by this
subgroup. This is a hyperbolic 3-orbifold with non-empty
singular locus. When $O$ has finite volume,
Theorem 1.1 gives a tower of finite-sheeted
covering spaces $\{ O_i \}$, and the subgroups of
$\pi_1(O)$ corresponding to these covers are
the required nested sequence. Suppose now that $O$
has infinite volume. Selberg's Lemma, together
with the assumption that $\pi_1(O)$ is finitely generated,
implies that $\pi_1(O)$
has a finite-index normal subgroup that is torsion free.
It is therefore isomorphic
to the fundamental group of a compact orientable
irreducible 3-manifold $M$ with non-empty boundary.
By the result of Cooper, Long and Reid (Theorem 1.3 of [4]),
$\pi_1(M)$ is large, unless $M$ is finitely covered
by an $I$-bundle over a surface with non-negative
Euler characteristic. But, $\pi_1(O)$ is then
elementary, contrary to assumption. Thus, as
argued in the bounded case in the proof of
Theorem 1.1, this implies the existence of
the required sequence of finite index subgroups.
$\square$

We are now in a position to prove Theorem 1.3
in the hyperbolic case. 

\noindent {\bf Theorem 1.3.} {\sl Let $O$ be a compact orientable 
geometric 3-orbifold
with non-empty singular locus. Then, the subgroup growth of $\pi_1(O)$ is
$$\cases{
\hbox{polynomial}, & if $O$ admits an $S^3$, ${\Bbb E}^3$, $S^2 \times {\Bbb E}$,
Nil or Sol geometry; \cr
\hbox{at least exponential}, & otherwise.}$$}

\noindent {\sl Proof (hyperbolic case).}
Let $O$ be a compact orientable 3-orbifold that has
non-empty singular locus and that is hyperbolic. Recall that this
means that the interior of $O$ admits a complete hyperbolic
structure, which may have finite or infinite volume. If $\pi_1(O)$ is
elementary, then $O$ also admits a Euclidean
structure and $\pi_1(O)$ has polynomial subgroup
growth. (See Section 8 for this deduction.)
If $\pi_1(O)$ is non-elementary and
$O$ has infinite volume, then, as in the proof
of Theorem 1.2, $\pi_1(O)$ is large, and so has (faster than)
exponential subgroup growth.

Thus, we may assume that $O$ has finite volume.
By Theorem 1.1,
there is an infinite nested sequence of finite-sheeted covers
$\{ O_i \}$ where $\pi_1(O_i)$ has linear growth
of mod $p$ homology for some prime $p$, and where
the degree of $O_{i+1} \rightarrow O_i$ is $p$
for each $i$. The fact that
$\pi_1(O)$ has at least exponential subgroup growth
then follows from the following proposition. $\square$

\noindent {\bf Proposition 5.1.} {\sl Let $G$ be a finitely generated
group, and let $p$ be a prime. Suppose that $G$ has an infinite 
nested sequence $\{ G_i \}$ of finite index subgroups 
with linear growth of mod $p$ homology. Suppose also
that the index $[G_i:G_{i+1}]$ is bounded above
independent of $i$. Then $G$ has
at least exponential subgroup growth.
Furthermore, if each $G_i$ is subnormal in $G$,
then $G$ has exponential subnormal subgroup growth.}

\noindent {\sl Proof.}  The number of
normal subgroups of $G_i$ with index $p$ is
$(p^{d_p(G_i)} -1)/(p-1)$. This is
a lower bound for the number of subgroups of $G$ with
index $p [G:G_i]$. Adjoining $G_i$ into this count, we deduce that
for $n = p [G:G_i]$,
$$s_n(G) \geq {p^{d_p(G_i)} -1 \over p-1} + 1> p^{d_p(G_i)-1}.$$
We are assuming that
$\inf_i d_p(G_i) / [G:G_i]$ is
some positive number $\lambda$. Hence,
$s_n(G) \geq p^{\lambda [G:G_i]-1}$.

We need to find a lower bound on $s_n(G)$ for
arbitrary positive integers $n$. This is
where we use the assumption that $[G_i:G_{i+1}]$ is bounded
above by some constant $k$. Thus, if we let
$[[n]]$ denote the largest integer less
than or equal to $n$ of the form $p [G:G_i]$,
we have the inequality $[[n]] \geq n/k$.
Therefore,
$$\liminf_n {\log s_n(G) \over n} 
\geq \liminf_n {\log s_{[[n]]}(G) \over k[[n]]}
\geq {\lambda \log p \over kp} > 0.$$
Thus, $G$ has at least exponential subgroup growth. 
In the case where each $G_i$ is subnormal in $G$,
the subgroups we are counting here are also subnormal,
and so $G$ then has exponential subnormal subgroup growth.
$\square$

\noindent {\bf Remark.} One might be tempted to think that Theorem 1.3
may be strengthened, because its proof appears to be quite wasteful
in the way it counts subgroups. In particular, it uses the existence
of only one tower of finite covers $O_i \rightarrow O''$, but in
fact many such towers are known to exists, that have the property that
the number of components of the singular locus of $O_i$ grows linearly in the degree of the cover.
This is because these covers were constructed from a nested
sequence of finite index subgroups of the group $\Gamma$.
Such a sequence was proved to exist by Theorem 4.4, using the fact that
$\Gamma$ fails the Golod-Shafarevich condition. However, it
is known that if a group fails this condition, then it
has many finite index subgroups (see Theorem 4.6.4 in [16] for example).
However, this does not lead to any significant improvement in
the subgroup growth of $\pi_1(O)$, since the known lower
bounds for the subgroup growth of $\Gamma$ are swamped
by the exponential terms arising from the linear growth
of mod $p$ homology of $\{ \pi_1(O_i) \}$.

We would now like to establish the following result, which is
a stronger version of Theorem 1.3. We will first prove it in the
case where the orbifold admits a finite volume hyperbolic structure. The
remaining seven geometries will be dealt with in the final
chapter.

\noindent {\bf Theorem 5.2.} {\sl Let $O$ be a compact orientable 
geometric 3-orbifold
with non-empty singular locus. Then, the subnormal subgroup growth of $\pi_1(O)$ is
$$\cases{
\hbox{polynomial}, & if $O$ admits an $S^3$, ${\Bbb E}^3$, $S^2 \times {\Bbb E}$,
Nil or Sol geometry; \cr
\hbox{exponential}, & otherwise.}$$}

We will need the following lemma.

\noindent {\bf Lemma 5.3.} {\sl Let $\{ G_i \}$ be a sequence of
finite index subgroups of a finitely generated group $G$,
and let $H$ be a finite index subnormal subgroup of $G$. If $\{ G_i \}$ has linear
growth of mod $p$ homology for some prime $p$, then $\{ G_i \cap H \}$
does also, after possibly discarding finitely many subgroups $G_i \cap H$.}

\noindent {\sl Proof.} An obvious induction allows 
us to reduce to the case where $H$ is normal in $G$.
Now, $G_i / (G_i \cap H)$ is isomorphic to $G_i H / H$,
which is a subgroup of $G/H$.
This places a uniform upper bound, independent of $i$,
on $[G_i : (G_i \cap H)]$ and, since there
are only finitely many subgroups of $G/H$, a
uniform upper bound on $d_p(G_i / (G_i \cap H))$.
Hence, by Proposition 3.6, there is a uniform
upper bound on $d_p(G_i) - d_p(G_i \cap H)$.
Thus, $\liminf_i d_p(G_i \cap H)/[G:G_i \cap H]$ is positive.
Hence the infimum of $d_p(G_i \cap H)/[G:G_i \cap H]$ is also
positive, once some initial terms of the sequence have been
deleted. The lemma now follows. $\square$

\noindent {\sl Proof of Theorem 5.2 (finite volume hyperbolic case).} 
When $O$ has boundary, the covering spaces
$\{ O_i \}$ we constructed in the proof of Theorem 1.1
had fundamental groups that were subnormal in $\pi_1(O)$,
and successive covers $O_{i+1} \rightarrow O_i$
had degree $p$. Hence by Proposition 5.1, $\pi_1(O)$
has exponential subnormal subgroup growth.

Suppose now that $O$ is closed.
In the proof of Theorem 1.1, we considered a tower
of covers $O_i \rightarrow O'' \rightarrow O' \rightarrow O$.
But we do not necessarily know that $O'' \rightarrow O'$
and $O' \rightarrow O$ are regular, and so we do not
know that $\pi_1(O_i)$ is a subnormal subgroup of
$\pi_1(O)$. Thus Proposition 5.1 cannot be applied
directly to deduce that $\pi_1(O)$ has exponential
subnormal subgroup growth.

Now, by Proposition 2.1, there is a finite regular cover
$M \rightarrow O'$ such that the composite cover
$M \rightarrow O' \rightarrow O$ is regular.
Hence, $\pi_1(M)$ is normal in $\pi_1(O')$ and
$\pi_1(O)$. According to Theorem 3.1, there is
a finite regular cover $\tilde O' \rightarrow O''$ such that
the composite cover $\tilde O' \rightarrow O'' \rightarrow O'$
is regular. Hence, $\pi_1(\tilde O')$ is normal in
$\pi_1(O')$ and $\pi_1(O'')$. Therefore, $\pi_1(M)
\cap \pi_1(\tilde O')$ is normal in $\pi_1(O')$, $\pi_1(\tilde O')$
and $\pi_1(M)$. Now, $\pi_1(O_i)$
is subnormal in $\pi_1(O'')$ and so $\pi_1(\tilde O')
\cap \pi_1(O_i)$ is subnormal in $\pi_1(\tilde O')$.
This implies that $\pi_1(M) \cap \pi_1(\tilde O') \cap
\pi_1(O_i)$ is subnormal in $\pi_1(M) \cap \pi_1(\tilde O')$.
We therefore have the chain of subgroups
$\pi_1(M) \cap \pi_1(\tilde O') \cap \pi_1(O_i)
\triangleleft \triangleleft \pi_1(M) \cap \pi_1(\tilde O')
\triangleleft \pi_1(M) \triangleleft \pi_1(O)$.
So each element of the sequence $\{
\pi_1(M) \cap \pi_1(\tilde O') \cap \pi_1(O_i) \}$ is
subnormal in $\pi_1(O)$. 

We claim that this sequence has linear growth of
mod $p$ homology, after one has possibly discarded some
initial terms in the sequence. Now, by Theorem 3.1, $\pi_1(\tilde O')$
is normal in $\pi_1(O'')$ and has finite index.
By Proposition 2.1, $\pi_1(M)$ is normal in
$\pi_1(O')$ and has finite index. Hence, $\pi_1(M) \cap
\pi_1(\tilde O')$ is normal in $\pi_1(\tilde O')$
and has finite index. Therefore, $\pi_1(M) \cap \pi_1(\tilde O')$
is subnormal in $\pi_1(O'')$ and has finite index.
Applying Lemma 5.3
(with $G = \pi_1(O'')$ and $G_i = \pi_1(O_i)$
and $H  = \pi_1(M) \cap \pi_1(\tilde O')$), we deduce
that this sequence has linear growth of mod $p$ homology, after
possibly discarding some initial terms of the sequence.
The theorem now follows by Proposition 5.1. $\square$

\vskip 18pt
\centerline{\caps 6. Property $(\tau)$, large groups and linear
growth of homology}
\vskip 6pt

In this section, we investigate to what extent the following
recent result [8] can be used to establish that the fundamental
group of a closed hyperbolic 3-orbifold $O$ with non-empty singular
locus is large.

\noindent {\bf Theorem 1.7.} {\sl Let $G$ be a finitely presented group,
let $p$ be a prime and suppose that $G \geq G_1 \triangleright G_2 \triangleright \dots$
is a nested sequence of finite index subgroups, such
that $G_{i+1}$ is normal in $G_{i}$ and has index a power of $p$, for each $i$. 
Suppose that $\{ G_i \}$ has linear growth of mod $p$ homology. Then,
at least one of the following must hold:q
\item{(i)} $G$ is large;
\item{(ii)} $G$ has Property $(\tau)$ with respect to $\{G_i \}$.

}

Let $O''$ be as in the proof of Theorem 1.1, and let
$G$ be $\pi_1(O'')$. As in Theorem 4.1,
let $C$ be a component of ${\rm sing}^c_p(O'')$, and 
let $\langle \! \langle \pi_1(\partial N(C))
\rangle \! \rangle$ be the subgroup of $\pi_1(O'')$ normally
generated by $\pi_1(\partial N(C))$. 
Let $\Gamma$ be $\pi_1(O'')/\langle \! \langle \pi_1(\partial N(C))
\rangle \! \rangle$.
It is shown in the proof of
Theorem 4.1 that $\Gamma$ has a presentation where
the inequality of the GS-$\tau$ Conjecture is satisfied.
Let us suppose that this conjecture is true.
It would imply that
$\Gamma$ does not have Property $(\tau)$
with respect to some nested sequence 
$\{ \Gamma_i \}$ of normal subgroups, each with
index a power of $p$. Let $G_i$ be the inverse image
of  $\Gamma_i$ in $G$.
Then, $G$ does not have Property $(\tau)$
with respect to $\{ G_i \}$.  As in the proof
of Theorem 4.1, let $O_i$
be the covering space of $O''$ corresponding
to $G_i$. It is shown there that
$|{\rm sing}^c_p(O_i)|$ is at least the degree 
of the cover $O_i \rightarrow O''$. By Proposition 4.3,
$\{ G_i \}$ therefore has linear growth of
mod $p$ homology. So, Theorem 1.7 implies that $G$ is large, which implies
that $\pi_1(O)$ is large. Thus, we
have proved the following.

\noindent {\bf Theorem 1.9.} {\sl The GS-$\tau$
Conjecture implies that the fundamental group of 
every closed hyperbolic 3-orbifold with non-empty singular locus
is large.}

\vskip 18pt
\centerline{\caps 7. Manifolds commensurable with hyperbolic 3-orbifolds}
\vskip 6pt

In this section, we consider manifolds $M$ that are commensurable with
a hyperbolic 3-orbifold $O$ with non-empty singular locus.
The aim is prove that many of the properties we have
deduced for $O$ also hold for $M$. In particular,
our goal is to prove Theorem 1.14.

\noindent {\bf Lemma 7.1.} {\sl Let $G$ be a finitely
generated group with exponential subnormal subgroup growth.
Then any finite index normal subgroup $H$ of $G$ also
has exponential subnormal subgroup growth.}

\noindent {\sl Proof.} If $G_i$ is a finite index subnormal subgroup
of $G$, then $G_i \cap H$ is subnormal in $H$. The index
$[H : G_i \cap H]$ is at most $[G:G_i]$.
For any given
subgroup of $H$ with index $n$, the number of ways of
writing it as $G_i \cap H$ for some subgroup
$G_i$ of $G$ is at most $(mn)^{\log m}$, where
$m = [G:H]$
(see the proof of Corollary 1.2.4 in [16]). Hence,
the fact that $G$ has exponential
subnormal subgroup growth implies that
$H$ does also. $\square$

\noindent {\bf Theorem 1.14.} {\sl Let $M$ be a 3-manifold
that is commensurable with a compact orientable finite-volume hyperbolic 3-orbifold $O$
with non-empty singular locus. Let $p$ be a prime that
divides the order of a torsion element of $\pi_1(O)$. Then $\pi_1(M)$ has
exponential subnormal subgroup growth, and $M$ has a nested
sequence of finite-sheeted covers that have linear growth
of mod $p$ homology and have Property $(\tau)$.}

\noindent {\sl Proof.} Let $M'$ be the common finite cover of $O$ and $M$. 
We may find a finite cover $M''$ of $M'$ such that $M'' \rightarrow M'$
and $M'' \rightarrow O$ are both regular covers. We may find
a further finite cover $M'''$ of $M''$ such that
$M''' \rightarrow M''$ and $M''' \rightarrow M$
are regular.

By Theorem 5.2, $\pi_1(O)$ has 
exponential subnormal subgroup growth.
Lemma 7.1, applied twice, implies that $\pi_1(M''')$ does also,
and hence so does $\pi_1(M)$.

By Theorem 1.1, $O$ has a nested
sequence of finite-sheeted covers $\{ O_i \}$
such that $\{ \pi_1(O_i) \}$ has linear growth
of mod $p$ homology, and where
each $\pi_1(O_i)$ is normal in $\pi_1(O_{i-1})$ and
has index $p$. Let $M_i$
be the covering space of $O$ corresponding to the
subgroup $\pi_1(M''') \cap \pi_1(O_i)$. 
These cover $M'''$ which covers $M$.
By Lemma 5.3 (with $G =\pi_1(O)$, $G_i = \pi_1(O_i)$ and
$H = \pi_1(M''')$), these covers have linear growth
of mod $p$ homology, after possibly discarding some
initial terms in the sequence.

Now, each $\pi_1(M_i)$ is subnormal in $\pi_1(M_1)$ and
has index a power of $p$. Hence, by Theorem 1.12
(with $G = \pi_1(M_1)$),
we may arrange that, in addition, this sequence of
subgroups has Property $(\tau)$.
$\square$

\vskip 18pt
\centerline{\caps 8. Geometric non-hyperbolic 3-orbifolds}
\vskip 6pt

In this section, we study compact orientable 3-orbifolds $O$ that admit a geometric
structure other than hyperbolic. Our goal is to prove Theorem 5.2 (and hence
Theorem 1.3) in this case. 

Note first that any compact geometric
3-orbifold is {\sl very good}: it admits a finite-sheeted
manifold cover (Corollary 2.27 of [3]). We start by considering orbifolds that admit a
geometry based on ${\Bbb H}^2 \times {\Bbb E}$ or
$\widetilde{PSL_2({\Bbb R})}$, but which do not admit
any of the remaining 6 geometries. Pass to a finite-sheeted
manifold cover $M$. Then,
$M$ is Seifert fibred, and the base orbifold is
hyperbolic. The base orbifold therefore admits
a finite-sheeted cover that is an orientable surface $S$ with
negative Euler characteristic. This pulls back to give
a finite covering $\tilde M \rightarrow M$.
Now, the Seifert fibration induces a surjective homomorphism
$\pi_1(\tilde M) \rightarrow \pi_1(S)$, and $\pi_1(S)$
admits a surjective homomorphism onto ${\Bbb Z} \ast {\Bbb Z}$.
Therefore, $\pi_1(O)$ is large. In particular,
its subgroup growth and subnormal subgroup growth
are (at least) exponential.

Suppose now that the orbifold admits a geometry based on 
$S^3$, ${\Bbb E}^3$, $S^2 \times {\Bbb E}$,
Nil or Sol geometry.
Any 3-orbifold modelled on spherical geometry
clearly has finite fundamental group, and hence
its subgroup growth is trivially polynomial.
When the model geometry is ${\Bbb E}^3$, Nil or Sol,
the orbifold is finitely covered by a torus bundle over
the circle. Hence, its fundamental group $G$ has a chain
of subgroups
$$G \triangleright G_0 \triangleright G_1 \triangleright \dots 
\triangleright G_n = \{ e \} \eqno (1)$$
where $|G/G_0|$ is finite and $G_i/G_{i+1} \cong {\Bbb Z}$.
It is a well-known fact, which is easy to prove
(see Corollary 1.4.3 of [16] for example), that this implies that $G$
has polynomial subgroup growth. The same argument
allows us to deal with orbifolds modelled on
$S^2 \times {\Bbb E}$ geometry. Thus, this proves
Theorem 5.2 for compact orientable 3-orbifolds
that admit a non-hyperbolic geometry.

In the cases above where the subgroup growth of the
fundamental group $G$ is polynomial, it is natural
to enquire about the {\sl degree} of this growth.
This is defined to be
$$\alpha(G) = \limsup_n {\log s_n(G) \over \log n}.$$
The determination of this quantity is likely to be
a tractable problem, but it is by no means trivial.
We merely mention here a few remarks about it.

It can be shown (Proposition 5.6.5 of [16]) that $\alpha(G) \leq h(G) + 1$,
where $h(G)$ is the Hirsch length of $G$. This
is defined to be the integer $n$ in the sequence (1) above.
In all the cases $G$ we considered, $h(G)$ is at most 3,
and therefore $\alpha(G)$ is at most 4.

Our second note is that $\alpha(G)$ is not necessarily
unchanged on passing to a finite index subgroup $H$.
It can be shown (Proposition 5.6.4 of [16]) that $\alpha(H)$
lies between $\alpha(G)$ and $\alpha(G) + 1$. It is 
slightly surprising to note that the upper bound is
sometimes realised. For example, when $G$ is the infinite
dihedral group, $\alpha(G) = 2$, but $G$ contains ${\Bbb Z}$
as an index two normal subgroup, and $\alpha({\Bbb Z}) = 1$.

It is clear that the degree of polynomial subgroup
growth for these 3-orbifold groups is worthy of
further investigation.

\vskip 18pt
\centerline{\caps References}
\vskip 6pt

\item{1.} {\caps M. Boileau, B. Leeb, J. Porti},
{\sl Geometrisation of 3-dimensional orbifolds},
Ann. Math. 162 (2005) 195--290.

\item{2.} {\caps J. Button},
{\sl Strong Tits alternatives for compact 3-manifolds with boundary,}
J. Pure Appl. Algebra 191 (2004) 89--98.

\item{3.} {\caps D. Cooper, C. Hodgson, S. Kerckhoff},
{\sl Three-dimensional Orbifolds and Cone-Manifolds},
MSJ Memoirs, Vol. 5 (2000)

\vfill\eject
\item{4.} {\caps D. Cooper, D. Long, A. Reid,} 
{\sl Essential closed surfaces in bounded $3$-manifolds,} 
J. Amer. Math. Soc. 10 (1997) 553--563.

\item{5.} {\caps D. Gabai, G.R. Meyerhoff, N. Thurston},
{\sl Homotopy hyperbolic 3-manifolds are hyperbolic}, 
Ann. Math. 157 (2003) 335--431.

\item{6.} {\caps M. Lackenby}, {\sl Heegaard splittings,
the virtually Haken conjecture and Property ($\tau$)},
Invent. Math. 164 (2006) 317--359.

\item{7.} {\caps M. Lackenby}, {\sl Some 3-manifolds and 3-orbifolds
with large fundamental group}, Preprint.

\item{8.} {\caps M. Lackenby}, {\sl Large groups, Property $(\tau)$ and
the homology growth of subgroups}, Preprint.

\item{9.} {\caps M. Lackenby, D. Long, A. Reid}, {\sl Covering spaces
of arithmetic 3-orbifolds,} Preprint.

\item{10.} {\caps A. Lubotzky}, {\sl Group presentations, $p$-adic analytic groups 
and lattices in SL (2,C),} Ann. Math. 118 (1983) 115--130.

\item{11.} {\caps A. Lubotzky}, {\sl Dimension function for discrete groups,} 
Proceedings of groups -- St. Andrews 1985,  254--262, London Math. Soc. Lecture Note Ser., 
121, Cambridge Univ. Press, Cambridge (1986).

\item{12.} {\caps A. Lubotzky}, {\sl Discrete Groups,
Expanding Graphs and Invariant Measures}, Progress
in Math. {\bf 125} (1994)

\item{13.} {\caps A. Lubotzky}, {\sl Eigenvalues of the 
Laplacian, the first Betti number and the congruence subgroup
problem,} Ann. Math. {\bf 144} (1996) 441--452.

\item{14.} {\caps A. Lubotzky}, Private communication.

\item{15.} {\caps A. Lubotzky, A. Mann,}
{\sl Powerful $p$-groups II. $p$-adic analytic groups,}
J. Algebra, 105 (1987) 506--515.

\item{16.} {\caps A. Lubotzky, D. Segal},
{\sl Subgroup growth}. Progress in Mathematics, 212. 
Birkh\"auser Verlag (2003)

\item{17.} {\caps J. Ratcliffe}, {\sl Euler characteristics of
3-manifold groups and discrete subgroups of ${\rm SL}(2, {\Bbb C})$},
J. Pure Appl. Algebra 44 (1987) 303--314.

\item{18.} {\caps P. Scott,} {\sl The geometries of $3$-manifolds}. 
Bull. London Math. Soc. 15 (1983) 401--487.

\item{19.} {\caps W. Thurston}, {\sl The geometry and topology
of 3-manifolds}, Lecture Notes, Princeton, 1978.

\vskip 6pt
\+ Mathematical Institute, University of Oxford, \cr
\+ 24-29 St Giles', Oxford OX1 3LB, UK. \cr

\end